\documentclass[11pt]{article}
\usepackage{setspace}
\usepackage{amsfonts}
\usepackage{amssymb}         
\usepackage{amsmath}    
\usepackage{multirow}
\usepackage{graphicx}
\usepackage{algorithm}
\usepackage{algorithmic}

\newtheorem{theorem}{Theorem}
\newtheorem{lemma}{Lemma}

\newtheorem{proposition}{Proposition}
\newtheorem{remark}{Remark}
\newtheorem{corollary}{Corollary}

\newcommand{\Pb}{{\mathbf{Pr}}}

\newcommand{\Pcal}{\mathcal{P}}

\begin{document}

\title{On the Exponential Probability Bounds for the Bernoulli Random Variables}
\author{\textsc{Vladimir Nikulin\thanks{Email: vnikulin.uq@gmail.com}}    \\ 
Department of Mathematics, University of Queensland \\
Brisbane, Australia
}
\date{}
\maketitle
\thispagestyle{empty}

\begin{abstract}
   We consider upper exponential bounds for the probability
   of the event that an absolute deviation of sample mean
   from mathematical expectation $p$ is bigger comparing with some ordered
   level $\varepsilon$.
   These bounds include 2 coefficients $\{\alpha, \beta\}$. 
In order to optimize the bound we are
   interested to minimize linear coefficient $\alpha$ and to 
maximize exponential coefficient $\beta$.
   Generally, the value of linear coefficient $\alpha$ may not be smaller than one.
   The following 2 settings were proved:
    1) $\{1, 2\}$ in the case of classical discreet problem as it was
   formulated by Bernoulli in the 17th century, and 2) 
$\{1, \frac{2}{1+\varepsilon^2}\}$ in the
   general discreet case with arbitrary rational $p$ and $\varepsilon.$
   The second setting represents a new structure of the exponential bound which may be
   extended to continuous case.
\end{abstract}

\section{Introduction}
\label{intro}

\noindent
Let $X, X_1, X_2, \cdots$ be a sequence of independent and identically distributed random
variables  $\Pb(X=1)=p, \Pb(X=0)=1-p.$

Ya. Bernoulli proved \cite{Bern86} that
$$\Pb(|\frac{1}{n} \sum_{i=1}^n X_i - p| > \varepsilon) \leq \frac{1}{1+C},$$
where
$$C=\min{(\frac{1}{s-1}\exp{( \xi_1 \log{\frac{r+1}{r}})},
     \frac{1}{r-1}\exp{( \xi_2 \log{\frac{s+1}{s}}))}},$$
$$\xi_1=\frac{k(r+1)+s}{r+s+1}, \xi_2=\frac{k(s+1)+r}{r+s+1},$$
$p=\frac{r}{r+s}, n=k(r+s), \varepsilon=\frac{1}{r+s},$  and
$k, r, s$ are arbitrary natural numbers.

As far as $\xi_1, \xi_2 \geq \frac{k}{r+s} = n \varepsilon^2$ we can make conclusion that
$$C=\min{(\frac{1}{s-1}\exp{( n \varepsilon^2 \log{\frac{r+1}{r})}},
     \frac{1}{r-1}\exp{( n \varepsilon^2 \log{\frac{s+1}{s}}))}}.$$

Therefore, we can formulate the law of large numbers
$$\lim_{n \rightarrow \infty} \Pb(|\frac{1}{n} 
\sum_{i=1}^n X_i - p| > \varepsilon) = 0,$$
where $r$ and $s$ are assumed to be fixed and $k \rightarrow \infty.$

Markov \cite{Mark14}, \cite{Mark24} considered case of arbitrary
$n, p$ and $\varepsilon$.
Uspensky \cite{Usp37} extended results of Markov further and derived the first
exponential bound
\begin{equation} \label{eq:uspen}
  \Pb(|\frac{1}{n} \sum_{i=1}^n X_i - p| \geq \varepsilon) \leq
  \alpha \exp{\{-\beta \varepsilon^2 n\}}
\end{equation}
with coefficients $\alpha = 2, \beta = 0.5.$
Additional and more detailed historical notes may be found in the Section
``Existing Inequalities'' \cite{Ben62}.

Hoeffding \cite{Hoef63} developed methods of \cite{Usp37}, \cite{Chern52} and
proved generally that $\alpha=\beta=2.$
Note that similar exponential bounds for the empirical distribution functions
were presented in \cite{Dvor56} and \cite{Mass90}.

In the Section~\ref{sec:bernoulli}
we prove that the following values of the constants $\alpha=1$ and $\beta=2$
may be used in the bound (\ref{eq:uspen})
in the discreet case formulated and considered by Ya. Bernoulli.
It is demonstrated in the Section~\ref{ssec:optval} that value of $\alpha$ can not be
 smaller than one. The following Section~\ref{ssec:oneside} proves
one-sided inequalities using methods and results of the Theorem~\ref{th:teorem}.

Section~\ref{sec:discreet} introduces a new structure of the exponential bound in the
general discreet case with arbitrary rational $p$ and $\varepsilon.$
The best bound in asymptotical sense corresponds to the bigger value of the
exponential coefficient. However, as it is discussed in the Section~\ref{ssec:relbetbnds},
the value of linear coefficient is also very important because of the practical reasons.

Section~\ref{sec:asymptot} represents an extension of the methods
and results of the Section~\ref{sec:discreet} to the continuous case.

Section~\ref{ssec:normal} considers particular application of the bounds to the
normalized sum of random variables.
It is demonstrated that all 3 types of bounds are asymptotically equivalent.
In this way, advantage of the propose bounds comparing with
Hoeffding's bound is absolute.

Note that using results of \cite{Ben02} and \cite{Ben04} we can
extend bounds for Bernoulli random variables to the case of
arbitrary bounded random variables.

Also, we note paper \cite{LeoPer04} where similar exponential
bounds were constructed for Markov chains. Generally, exponential
bounds proved to be very effective in order to define required
size of the sample in order to ensure proper quality of
estimation, see, for example,
\cite{Lue98,HesZat00,GamMed04,RavLaf04,YehGal05}.

\section{Main lemma and definitions} \label{sec:fsec}

We will use essentially different approach comparing with
 \cite{Hoef63} and  \cite{Chern52}. This technique is based on the properties
of convex (concave) functions applied to the binomial coefficients.
The following Lemma formulates the core of the methods
which are employed in  the Sections~\ref{sec:bernoulli},
~\ref{sec:discreet} and ~\ref{sec:asymptot}.

\begin{lemma} \label{lm:main}
Suppose that $\varphi(t), 1 \leq t \leq n,$ is a function with continuous
second derivative where $n$ is any natural number. Then,
\begin{equation} \label{eq:convf}
  \sum_{j=1}^n \varphi(j) \geq n \varphi(\frac{n+1}{2})
\end{equation}
if $\varphi'(t) \leq 0, \varphi''(t) \geq 0, 1 \leq t \leq n$ (convex case);
\begin{equation} \label{eq:concf}
  \sum_{j=1}^n \varphi(j) \geq n \frac{\varphi(1)+\varphi(n)}{2}
\end{equation}
if $\varphi'(t) \geq 0, \varphi''(t) \leq 0, 1 \leq t \leq n$ (concave case).
\end{lemma}
\textit{Proof:}
The following representations are valid
\begin{displaymath}   
  \sum_{j=1}^n \varphi(j) =  \left\{     
\begin{array}{ll}   
\sum_{j=0}^{m-1} \left[ \varphi(1+j) + \varphi(n-j) \right]
     \hspace{0.08in} if \hspace{0.05in} n=2m; \\
     \sum_{j=0}^{m-1} \left[ \varphi(1+j) + \varphi(n-j) \right] + \varphi(m+1)
     \hspace{0.08in} if \hspace{0.05in} n=2m+1,
\end{array} \right.
\end{displaymath}   
where $m$ is a natural number.
We obtain required bounds combining above equation with
$$\varphi(1+j)+\varphi(n-j) \geq 2 \varphi(\frac{n+1}{2}) 
\geq 2 \varphi(m+1), 1 \leq j \leq m-1,$$
in the convex case, and with
$$2\varphi(m+1) \geq \varphi(1+j)+\varphi(n-j) 
\geq \varphi(1)+\varphi(n), 1 \leq j \leq m-1,$$
in the concave case.
$\blacksquare$

The following notations will be used below
$$\Pcal^0 := \Pb(\frac{1}{n} \sum_{i=1}^n X_i = p) =
\left( \begin{array}{cc} n \\ m \end{array} \right) p^m (1-p)^{n-m},$$
$$\Pcal^{+}_1 := \Pb(0 < \frac{1}{n} \sum_{i=1}^n X_i - p \leq \varepsilon), \hspace{0.05in}
\Pcal^{+}_2 := \Pb(\varepsilon < \frac{1}{n} \sum_{i=1}^n X_i -p),$$
$$\Pcal^{-}_1 := \Pb(0 < p -  \frac{1}{n} \sum_{i=1}^n X_i \leq \varepsilon), \hspace{0.05in}
\Pcal^{-}_2 := \Pb(\varepsilon < p - \frac{1}{n} \sum_{i=1}^n X_i),$$
where $m=kr.$

Assuming that $n = k(r+s)$ we form groups of binomial probabilities of the equal size $k$
$$S_j :=\sum_{\ell=m-kj}^{m-1-(j-1)k}
\left( \begin{array}{cc} n \\ \ell \end{array} \right)
p^{\ell} (1-p)^{n-\ell}, j=1,\ldots,r \mbox{(left groups)},$$
$$Z_j := \sum_{\ell=m+1+(j-1)k}^{m+kj}
\left( \begin{array}{cc} n \\ \ell \end{array} \right)
 p^{\ell} (1-p)^{n-\ell}, j=1,\ldots,s \mbox{(right groups)}.$$

Then, we consider relations of the corresponding binomial coefficients from the
neighbor groups
\begin{subequations}
\begin{align}
  \label{eq:step2}
  A(j) := \frac{\left( \begin{array}{cc} n \\ m-j \end{array} \right)}
  {\left( \begin{array}{cc} n \\ m-k-j \end{array} \right)} =
  \prod_{v=1}^k \frac{n-m+j+v}{m-k-j+v}, j=1,\ldots,k(r-1), \\
  \label{eq:step1}
  B(j) := \frac{\left( \begin{array}{cc} n \\ m+j \end{array} \right)}
  {\left( \begin{array}{cc} n \\ m+k+j \end{array} \right)} =
  \prod_{v=1}^k \frac{m+j+v}{n-m-k-j+v}, j=1,\ldots,k(s-1),
\end{align}
\end{subequations}
where $s, r \geq 2.$

 By definition
\begin{equation} \label{eq:bnd2a}
  \Pb(|\frac{1}{n} \sum_{i=1}^n X_i - p| > \varepsilon) =
  \Pcal^+_2 + \Pcal^-_2
  = \sum_{j=2}^s Z_j +\sum_{j=2}^r S_j.
\end{equation}

\begin{remark} Note that $\Pcal^+_2=0$ if $s=1$, and $\Pcal^-_2=0$ if $r=1.$
\end{remark}

\section{Bernoulli problem}  \label{sec:bernoulli}

We exclude from consideration the trivial cases: $p=0$ and $p=1$,
and assume that $0 < p < 1.$
As it will be demonstrated below the task of estimation of $\Pcal^{+}_2$ is easier
 comparing with estimation of $\Pcal^{-}_2$ if $p > 0.5$. On the other hand,
the task of estimation of $\Pcal^{-}_2$ is easier comparing with estimation of $\Pcal^{+}_2$
if $p < 0.5$. As far as the problem is symmetrical, we assume that $p \geq 0.5.$

\begin{theorem} \label{th:teorem} Suppose that
 $p=\frac{r}{r+s}, n=k(r+s), \varepsilon=\frac{1}{r+s},$ and
$k, r, s$ are arbitrary natural numbers. Then,
$$b \cdot Z_{j+1} \leq Z_j, j=1,\ldots,s-1,$$
where $b := \exp{\{2 \varepsilon^2 n\}}$.
Therefore,
\begin{equation} \label{eq:bnd1}
  (b-1) \Pcal^{+}_2 \leq \Pcal^{+}_1.
\end{equation}
\end{theorem}
\textit{Proof:}
By (\ref{eq:step1}), $B(j)$ is an increasing function of $j \geq 1$.
Therefore,
$$\lambda_j := \frac{Z_j}{Z_{j+1}} \geq B((j-1)k+1)\left(\frac{1-p}{p}\right)^k=
\prod_{v=1}^k \frac{m+(j-1)k+1+v}{n-m-jk-1+v} \left(\frac{s}{r}\right)^k$$
or
$$\log{\lambda_j} \geq k\log{\frac{s}{r}} + 
\sum_{v=1}^k \varphi(j, v), j=1,\ldots,s-1,$$
where $$\varphi(j, v) = \log{\frac{m+(j-1)k+1+v}{n-m-jk-1+v}}$$
is a convex function of $v$ because $r \geq s$ or $m \geq n-m$, and
 according to (\ref{eq:convf})
$$\lambda_j \geq \left(\frac{s(m+jk-0.5k+1.5)}{r(n-m-jk+0.5k-0.5)}\right)^k$$
\begin{equation} \label{eq:bnd2}
 \geq \left( \frac{s(m+jk-0.5k+2)}{r(n-m-jk+0.5k)}\right)^k =(1+\delta)^k,
\end{equation}
where $$\delta=\frac{n(j-0.5)+2s}{kr(s-j+0.5)}.$$

Using general inequality
\begin{equation} \label{eq:gbound}
  \log{(1+\delta)} \geq \frac{2\delta}{2+\delta} \hspace{0.1in} \forall \delta \geq 0,
\end{equation}
we obtain
$$\log{\lambda_j} \geq \frac{2k(n(j-0.5)+2s)}{2(s-j+0.5)kr+n(j-0.5)+2s}$$
\begin{equation} \label{eq:bnd3}
 \geq \frac{2kn(j-0.5)}{2srk+(j-0.5)(n-2kr)}, j=1,\ldots,s-1.
\end{equation}

Furthermore, based on the following properties
$$2srk+(j-0.5)(n-2rk)=k\left[ r(2s-j+0.5)+s(j-0.5)\right]>0, \hspace{0.05in} 2kr \geq n,$$
we have
$$\log{\lambda_j} \geq \frac{n(j-0.5)}{sr}.$$
Then, using relations
$r=\frac{p}{\varepsilon}, s=\frac{1-p}{\varepsilon}$
we transform above inequality to the required form
\begin{equation} \label{eq:bnd4}
 \lambda_j \geq  \exp{\{\frac{n\varepsilon^2(2j-1)}{2p(1-p)}\}} \geq 
\exp{\{2n \varepsilon^2(2j-1)\}},
  j=1,\ldots,s-1.
\end{equation}
Above equation (\ref{eq:bnd4}) completes proof of the Theorem.
$\blacksquare$

\begin{theorem} \label{th:teorem2}
 Suppose that
 $p=\frac{r}{r+s}, n=k(r+s), \varepsilon=\frac{1}{r+s},$ and
$k, r, s$ are arbitrary natural numbers. Then,
$$b \cdot S_{j+1} \leq S_j, j=1,\ldots,r-1.$$
Therefore,
\begin{equation} \label{eq:bnd1a}
  (b-1) \Pcal^{-}_2 \leq \Pcal^{-}_1.
\end{equation}
\end{theorem}
\textit{Proof:}
By (\ref{eq:step2}), the coefficient $A(j)$ is an increasing function of $j \geq 1$.
Therefore,
\begin{equation} \label{eq:bnd5}
  q_j := \frac{S_j}{S_{j+1}} \geq A((j-1)k+1) \left( \frac{p}{1-p}\right)^k
\end{equation}
or
$$\log{q_j} \geq k \log{\frac{r}{s}} + \sum_{v=1}^k \psi(j, v),$$
where $\psi(j, v) = \log{\frac{n-m+(j-1)k+1+v}{m-jk-1+v}}.$

The following condition
$$k(r-s-1) \geq 2$$
give us a guarantee that the function $\psi(1, v)$ is concave.

Suppose that $r = s+2$ and $k \geq 2$ or $r \geq s+3$ and $k \geq 1.$
Then, by (\ref{eq:concf}) we have
$$\log{q_j} \geq \log{q_1} \geq k \log{\frac{r}{s}}+\frac{k}{2}
\log{\frac{(n-m+2)(n-m+k+1)}{(m-k)(m-1)}}$$
\begin{equation} \label{eq:bnd6a}
  \geq \frac{k}{2} \log{\frac{r(s+1)}{s(r-1)}} = \frac{k}{2}
\log{\left(1+\frac{r+s}{s(r-1)}\right)} \geq \frac{k(r+s)}{2s(r-1)+r+s}
\end{equation}
$$=\frac{n\varepsilon^2}{2p(1-p)+\varepsilon(2p-1)}\geq 2n\varepsilon^2$$
if $4p(1-p)+2\varepsilon(2p-1) \leq 1$ or $\varepsilon \leq \frac{1-4p(1-p)}{2(2p-1)}=p-0.5.$
Thus, $\varepsilon=\frac{1}{r+s} \leq \frac{r}{r+s}-0.5$ or $r \geq s+2.$

The case $r=s+2, k=1$ is easy to consider:
$$\log{q_1} = \log{\frac{(s+2)^2}{s(s+1)}} = \log{\left(1+\frac{3s+4}{s(s+1)}\right)} \geq
\frac{6s+8}{2s^2+5s+4}\geq \frac{\gamma}{2(s+1)} = \gamma n \varepsilon^2$$
or
$$(12-2\gamma)s^2+(28-5\gamma)s+16-4\gamma \geq 0.$$
Above inequality is valid $\forall s\geq 1$ if $\gamma \leq 4$.

Now, we consider remaining case $r = s+1$ (the case $r=s$ was considered already in
the first part of the proof because the problem is symmetrical).
Then, $\psi(j, v)$ is convex as a function of $v$, and, by (\ref{eq:convf}),

\begin{equation} \label{eq:bnd6}
  q_j \geq \left(\frac{(s+1)(ks+jk+1.5-0.5k)}{s(ks-jk-0.5+1.5k)}\right)^k =(1+\delta)^k,
\end{equation}
where $$\delta=\frac{2jsk+2s-ks+jk+1.5-0.5k}{s(ks-jk-0.5+1.5k)}.$$

Using inequality (\ref{eq:gbound}), we obtain
\begin{equation} \label{eq:bnd3}
  \log{q_j} \geq \frac{2k(2jsk+2s-sk+jk+1.5-0.5k)}{2s^2k+s+2sk+jk+1.5-0.5k} \geq
  \frac{\gamma k}{2s+1}=\gamma n \varepsilon^2
\end{equation}
or
$$(8jk+8-4k-2k\gamma)s^2+(8jk-4k+10-(2k+1)\gamma)s + 2jk + 3 - k$$
$$- (jk+1.5-0.5k)\gamma \geq 0.$$
Above inequality is valid for $\gamma \leq 2$ assuming that $j, s \geq 1.$
$\blacksquare$

\begin{corollary} \label{cor:corr}
The following upper bound is valid
\begin{equation} \label{eq:bnd1}
  \Pb(|\frac{1}{n} \sum_{i=1}^n X_i - p| > \varepsilon) <
  \exp{\{-2 \varepsilon^2 n\}},
\end{equation}
where $p=\frac{r}{r+s}, n=k(r+s), \varepsilon=\frac{1}{r+s},$ and
$k, r, s$ are arbitrary natural numbers.
\end{corollary}
\textit{Proof:}
It follows from the Theorems~\ref{th:teorem} and
\ref{th:teorem2}
that
\begin{equation} \label{eq:geom}
  \sum_{j=2}^r S_j \leq \frac{S_1}{b-1}, \hspace{0.07in}
  \sum_{j=2}^s Z_j \leq \frac{Z_1}{b-1}.
\end{equation}

The following representation is valid
$$\Pb(|\frac{1}{n} \sum_{i=1}^n X_i - p| > \varepsilon)=
\frac{\sum_{j=2}^r S_j + \sum_{j=2}^s Z_j}{\Pcal^0 +
\sum_{j=1}^r S_j + \sum_{j=1}^s Z_j}.$$

It easy to see that $f(x) =\frac{x}{C+x}$ is an increasing function of $x,$
where $C$ is a positive constant. Therefore,
$$\Pb(|\frac{1}{n} \sum_{i=1}^n X_i - p| > \varepsilon) <
\frac{\frac{S_1}{b-1}+\frac{Z_1}{b-1}}{S_1+Z_1+\frac{S_1}{b-1}+\frac{Z_1}{b-1}} =
\frac{1}{b}=\exp{(-2n\varepsilon^2)}.$$
$\blacksquare$

\begin{figure}[h]
\includegraphics[scale=0.75]{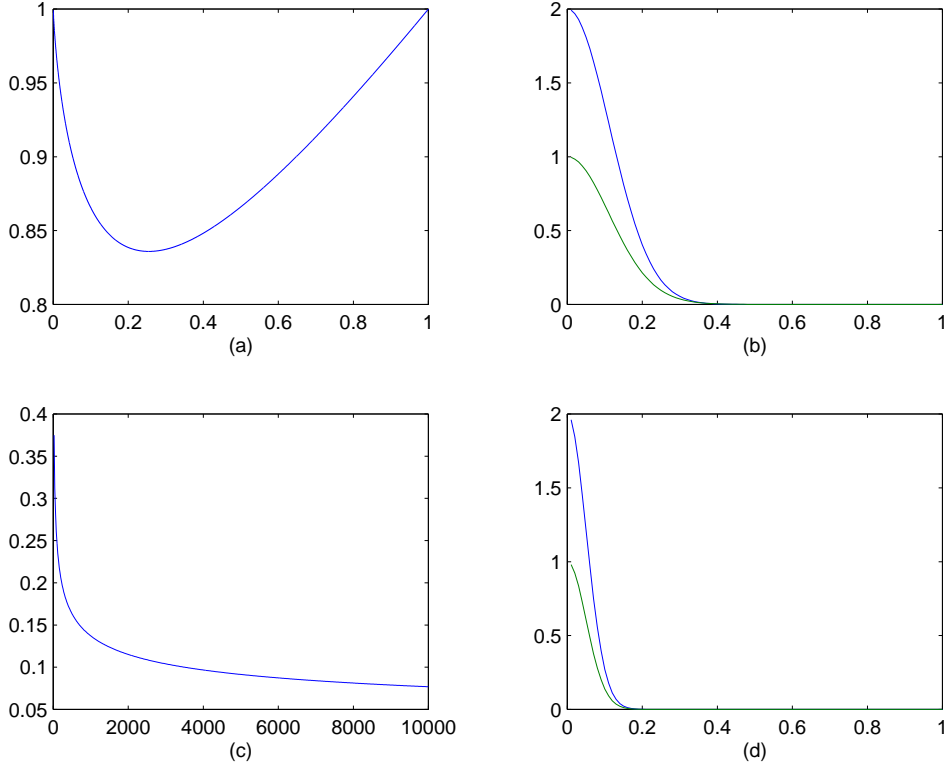}
  \caption{a) $f(p) = p^p (1+p)^{1-p}, 0 \leq p \leq 1;$
  c) function $\mu(n)$ which is defined in (\ref{eq:bnd63});
  bounds (\ref{eq:bnd56}) and (\ref{eq:hff5}) as a function of  $\varepsilon$
  where b) $n=20$, d) $n=100$.}
\label{fig: figure1}
\end{figure}

\subsection{Optimal value of the linear coefficient} \label{ssec:optval}

Clearly, value of the linear coefficient $\alpha=2$ in the exponential bound
may not be
regarded as an optimal. The following Proposition~\ref{pr:prop} demonstrates that
$\alpha \geq 1.$

\begin{proposition} \label{pr:prop}
Suppose that $0 < p < 1.$ Then, we can construct sequence of deviations
$\varepsilon(n) \underset{n \rightarrow \infty}{\longrightarrow} 0$ so that
$\varepsilon^2 n \underset{n \rightarrow \infty}{\longrightarrow} 0$ and
$\Pcal^-_2 + \Pcal^+_2 \underset{n \rightarrow \infty}{\longrightarrow} 1.$
\end{proposition}
\textit{Proof:}
By definition $\Pcal^0 = \left( \begin{array}{cc} n \\ m \end{array} \right)
p^m (1-p)^{n-m}.$
The following relation is valid
$$\log{\left( \begin{array}{cc} n \\ m \end{array} \right)} =
\sum_{j=1}^{n-m} \log{\frac{m+j}{j}}.$$
By Lemma~\ref{lm:main},
 $$\log{\left( \begin{array}{cc} n \\ m \end{array} \right)} \leq
(n-m)\log{\frac{1+n+m}{1+n-m}} \leq (n-m)\log{\frac{n+m}{n-m}} = (n-m)\log{\frac{1+p}{1-p}}.$$
Respectively,
$$\Pcal^0 \leq \left( p^p (1+p)^{1-p}\right)^n,$$
where $0< p^p (1+p)^{1-p} <1 \hspace{0.08in} \forall p:  0 < p < 1.$
Figure~\ref{fig: figure1}(a) illustrates graph of the function
$f(p) = p^p (1+p)^{1-p}, 0 \leq p \leq 1.$

Suppose that $\varepsilon = \frac{1}{2n}$. 
Then, $\Pcal^+_1$ and $\Pcal^-_1$ will be empty sets,
and $\Pcal^-_2 + \Pcal^+_2 \underset{n \rightarrow \infty}{\longrightarrow} 1$ as far as
$\Pcal^0 \underset{n \rightarrow \infty}{\longrightarrow} 0.$
At the same time, by definition,
\newline
$\varepsilon^2 n \underset{n \rightarrow \infty}{\longrightarrow} 0.$
$\blacksquare$

\begin{corollary} Coefficient $\alpha$ in the bound (\ref{eq:uspen}) 
may not be smaller than 1 for the arbitrary $n \geq 2$ and $\varepsilon > 0.$
\end{corollary}

\subsection{One-sided inequalities} \label{ssec:oneside}

In this Section we will use again the property that left groups of binomial 
probabilities may be
estimated more effectively if $p < 0.5.$ Analogously,
right groups of binomial probabilities may be estimated more effectively if $p > 0.5.$

\begin{proposition} \label{pr:prop5}
Suppose that $p = \frac{r}{r+s}, n=k(r+s), 
\delta = \nu \varepsilon = \frac{\nu}{r+s},$ 
where $k, r, s$ and $\nu$ are arbitrary natural numbers.
Then,
$$\Pb(\frac{1}{n} \sum_{i=1}^n X_i -p > \delta) \leq \frac{1}{2}
  \exp{( -\frac{n\delta^2}{2p(1-p)})}$$
if $r \geq s,$ ($p \geq 0.5$);

$$\Pb(p - \frac{1}{n} \sum_{i=1}^n X_i > \delta) \leq \frac{1}{2}
  \exp{(-\frac{n\delta^2}{2p(1-p)})}$$
if $r \leq s,$ ($p \leq 0.5$).
\end{proposition}
\textit{Proof:}
Combining the equality
$$\sum_{i=j}^{\nu} (2i-1) = \nu^2 - (j-1)^2, j=1,\ldots,\nu,$$
and (\ref{eq:bnd4}) we derive that
$$Z_j \geq Z_{\nu+1} \exp{\{\frac{n\varepsilon^2(\nu^2-(j-1)^2)}{2p(1-p)}\}}, j=1,\ldots,\nu,$$
and
\begin{equation} \label{eq:one1}
   \sum_{j=1}^{\nu} Z_j \geq Z_{\nu+1} \sum_{j=1}^{\nu}
  \exp{\{\frac{n\varepsilon^2(\nu^2-(j-1)^2)}{2p(1-p)}\}}.
\end{equation}
Next, we use the same bound (\ref{eq:bnd4}) in the opposite direction
$$Z_j \geq Z_{j+1} 
\exp{\{\frac{n\varepsilon^2(2\nu+1)}{2p(1-p)}\}}, j=\nu+1,\ldots,s-1,$$
and
\begin{equation} \label{eq:one2}
   \sum_{j=\nu+1}^{s} Z_j \leq Z_{\nu+1} \left( 1-
  \exp{\{-\frac{n\varepsilon^2(2\nu+1)}{2p(1-p)}\}} \right)^{-1}.
\end{equation}
Again, we make an assumption $r \geq s$ or $p \geq 0.5$
without loss of generality,
and construct upper bound for the conditional probability
$$\Pb(\frac{1}{n} \sum_{i=1}^n X_i -p > \delta) \mid
  \frac{1}{n} \sum_{i=1}^n X_i > p) = \frac{\sum_{j=\nu+1}^s Z_j}{ \sum_{j=1}^s Z_j}$$
\begin{equation} \label{eq:one3}
  \leq \left( \left(1-\exp{(-\frac{n\varepsilon^2(2\nu+1)}{2p(1-p)})}
  \right) \sum_{j=1}^{\nu} \exp{(\frac{n\varepsilon^2(\nu^2-(j-1)^2)}{2p(1-p)})}+1\right)^{-1},
\end{equation}
where the last inequality may be regarded as a consequence of (\ref{eq:one1}) 
and (\ref{eq:one2}).

Then, we simplify (\ref{eq:one3})
using properties $\nu^2-(j-1)^2 -2\nu-1 \leq \nu^2-j^2, j=1,\ldots,\nu,$
\begin{equation} \label{eq:one7}
  \Pb(\frac{1}{n} \sum_{i=1}^n X_i -p > \delta) \mid
  \frac{1}{n} \sum_{i=1}^n X_i > p) \leq  \exp{\{-\frac{n \delta^2}{2p(1-p)}\}}.
\end{equation}

As a next step we prove that
$$\Pb(\frac{1}{n} \sum_{i=1}^n X_i \leq p \mid p \geq 0.5) \geq \Pb(\frac{1}{n}
\sum_{i=1}^n X_i > p \mid p \geq 0.5).$$

Let us consider relations of the binomial probabilities which are symmetrical against
central point $m=kr$
$$\frac{\left( \begin{array}{cc} n \\ m - i +1 \end{array} \right)}
{\left( \begin{array}{cc} n \\ m+i \end{array}\right)} \left( \frac{1-p}{p}\right)^{2i-1}
 = \prod_{j=1}^{2i-1} \frac{(m-i+1+j)(n-m)}{(n-m-i+j)m}$$
$$\geq \left( \frac{m+1}{m} \right)^{2i-1} > 1, i=1,\ldots,n-m,$$
where we used (\ref{eq:convf}). Therefore,
\begin{equation} \label{eq:one8}
  \Pb(\frac{1}{n} \sum_{i=1}^n X_i \geq p) = \Pb(\frac{1}{n} 
\sum_{i=1}^n X_i \geq p \mid p \geq 0.5) \leq 0.5.
\end{equation}

By the way, using identical method, one can prove
$$\Pb(\frac{1}{n} \sum_{i=1}^n X_i \leq p \mid p \leq 0.5) \leq 0.5.$$

The final upper bound follows from (\ref{eq:one7}) and (\ref{eq:one8}) applied to the
Bayesian formulae
$$\Pb(\frac{1}{n} \sum_{i=1}^n X_i -p > \delta)$$
$$ = \Pb(\frac{1}{n} \sum_{i=1}^n X_i \geq p) \cdot
  \Pb(\frac{1}{n} \sum_{i=1}^n X_i -p > \delta) \mid
  \frac{1}{n} \sum_{i=1}^n X_i > p).$$
$\blacksquare$

\section{General discreet case} \label{sec:discreet}

In this section we assume that $p$ and $\varepsilon$
may be represented by positive
rational numbers with denominator as a number of observations in the sample.
Respectively, we will cover all possible empirical values of the sample 
mean as an estimator of the
probability $p$.
Note that the role of the parameter $k$ will be different here comparing 
with previous section.

\begin{theorem} \label{th:teorem3}
Suppose that
 $p=\frac{m}{n} \geq 0.5, \varepsilon=\frac{k}{n},$ and
$k, m, n$ are arbitrary natural numbers with condition $k, m < n.$ Then,
$$b \cdot Z_{j+1} \leq Z_j, j=1,\ldots,s-1,$$
where $b := \exp{\{2 \varepsilon^2 n\}}$.
Therefore,
\begin{equation} \label{eq:bnd2a}
  (b-1) \Pcal^{+}_2 \leq \Pcal^{+}_1.
\end{equation}
\end{theorem}
\textit{Proof:} 
The following upper bound may be obtained similarly to (\ref{eq:bnd2})
\begin{equation} \label{eq:bnd21}
  \log{\lambda_j} \geq k \log{\frac{(n-m)(m+0.5k+2)}{m(n-m-0.5k)}}, j=1,\ldots,s-1,
\end{equation}
where $2m \geq n$ and $p+\frac{\varepsilon}{2} < 1$ (else $\Pcal^+_2 = 0$).

It follows from ( \ref{eq:bnd21} ) that
\begin{equation} \label{eq:bnd22}
  \log{\lambda_j} \geq k \log{(1 + \delta)}
\end{equation}
where $\delta = \frac{0.5kn+2n-2m}{m(n-m-0.5k)}.$

By (\ref{eq:gbound}), we obtain
$$\frac{k(kn+4n-4m)}{2nm-2m^2-km+0.5kn+2n-2m} \geq \frac{\gamma k^2}{n} = 
\gamma \varepsilon^2 n$$
or
$$k+4-4p \geq \gamma k \left[ 2p(1-p)+\frac{\varepsilon}{2} (1-2p) \right] +
2 \gamma \varepsilon (1-p).$$

Above inequality is valid $\forall k \geq 1$ and
$\forall p \geq 0.5$ if $\gamma \leq 2.$
$\blacksquare$

The following Theorem~\ref{th:teorem4} introduces a new structure of the
exponential bound as a main result of this paper.

\begin{theorem} \label{th:teorem4} Suppose that
 $p=\frac{m}{n} \geq 0.5, \varepsilon=\frac{k}{n},$ and
$k, m, n$ are arbitrary natural numbers with conditions $2 \leq k <n, m < n.$ Then,
$$ \exp{\{\frac{2 \varepsilon^2 n}{1+\varepsilon^2}\}} 
\cdot S_{j+1} \leq S_{j}, j=1,\ldots,r-1.$$
Therefore,
\begin{equation} \label{eq:bth24}
  (\exp{\{\frac{2 \varepsilon^2 n}{1+\varepsilon^2}\}}-1) \Pcal^{-}_2 \leq \Pcal^{-}_1 
\end{equation}
(the case $k=1$ is covered by the Theorem~\ref{th:teorem2} with stronger result).
\end{theorem}
\textit{Proof:}
Similar to (\ref{eq:bnd6a}) and assuming that $2m > n+k +2$ we derive
\begin{equation} \label{eq:bth63}
  \log{q_j} \geq k \log{\frac{m}{n-m}}+ \frac{k}{2}
  \log{\frac{(n-m+2)(n-m+k+1)}{(m-k)(m-1)}}, 
\end{equation}
where $j=1,\ldots,r-1, m > k \geq 1$ (else $\Pcal^-_2 =0$).

The following relations are valid
$$\frac{(n-m+2)(n-m+k+1)}{(m-k)(m-1)} \left(\frac{m}{n-m}\right)^2 \geq
\frac{m(n-m+k+1)}{(m-k)(n-m)}=1+\delta,$$
where $\delta=\frac{m+kn}{(m-k)(n-m)}.$

By (\ref{eq:gbound}), we obtain
$$\log{q_j} \geq \frac{k(m+kn)}{2m(n-m+k)+m-kn}
\geq \frac{\gamma k^2}{n} = \gamma n \varepsilon^2$$
or
$$n(m+kn) \geq \gamma k \left[2m(n-m+k)+m-kn\right]$$
or
$$p+k \geq \gamma k \left[ 2p(1-p+\varepsilon) +\frac{p}{n} -\varepsilon\right]=
\gamma k \left[ 2p(1-p) + \varepsilon (2p-1) \right] + \gamma \varepsilon p.$$

Above inequality take place if
\begin{equation} \label{eq:bnd32}
  \gamma \leq \frac{2}{1+\varepsilon^2}.
\end{equation}

Let us consider remaining case $n \leq 2m \leq n+k +2.$
Similar to (\ref{eq:bnd2}) we have
\begin{equation} \label{eq:bnd31}
  q_j \geq \left(\frac{m(n-m+1.5+0.5k)}{(n-m)(m-0.5k-0.5)}\right)^k =(1+\delta)^k,
\end{equation}
where 
$$\delta=\frac{m+0.5kn+0.5n}{nm-m^2+0.5(km+m-kn-n)}.$$

By (\ref{eq:gbound}), we obtain
$$\frac{2k(m+0.5kn+0.5n)}{2nm-2m^2+km+2m-0.5n(k+1)} \geq \frac{\gamma k^2}{n} = 
\gamma \varepsilon^2 n$$
or
$$k+2p+1 \geq \gamma k (2p(1-p) +0.5 
\varepsilon (2p-1)) + \gamma \varepsilon (2p-0.5).$$
Above inequality take place if
\begin{equation} \label{eq:bnd33}
  \gamma \leq \frac{8}{4+\varepsilon^2}.
\end{equation}
Note that condition (\ref{eq:bnd33}) is less restrictive comparing with (\ref{eq:bnd32}).
$\blacksquare$

\begin{corollary} The following upper bound is valid
\begin{equation} \label{eq:bnd56}
  \Pb(|\frac{1}{n} \sum_{i=1}^n X_i - p| > \varepsilon) <
  \exp{\{-\frac{2 \varepsilon^2 n}{1+\varepsilon^2}\}},
\end{equation}
where $p=\frac{m}{n}, \varepsilon=\frac{k}{n},$ and
$k, m, n$ are arbitrary natural numbers with conditions $2 \leq k <n, n \leq 2m < 2n.$
\end{corollary}

Proof is similar to the Corollary~ \ref{cor:corr} : we obtain required value 
of the parameter
$\beta$ as $\min{\{2, \frac{2}{1+\varepsilon^2}, \frac{8}{4+\varepsilon^2}\}}$ 
where first value follows
from Theorem~\ref{th:teorem3}, second and third values follow from (\ref{eq:bnd32}) 
and (\ref{eq:bnd33}).

\subsection{Relations between bounds} \label{ssec:relbetbnds}

The bound (\ref{eq:bnd56}) improves Hoeffding's bound
\begin{equation} \label{eq:hff5}
  \Pb(|\frac{1}{n} \sum_{i=1}^n X_i - p| > \varepsilon) \leq 2 \exp{\{-2 \varepsilon^2 n\}}
\end{equation}
if
\begin{equation} \label{eq:bnd63}
  0 < \varepsilon^2 \leq \mu^2(n) = \frac{\phi}{2} +
   \sqrt{ \phi \left( 1+\frac{\phi}{4}\right)}, \hspace{0.08in} \phi := \frac{\log{2}}{2n}.
\end{equation}

Therefore,
$$\lim_{n \rightarrow \infty} n^{\frac{1}{4}} \mu(n) =
\left( \frac{\log{2}}{2}\right)^{\frac{1}{4}}.$$
Figure~\ref{fig: figure1}(c) illustrates behavior of the function $\mu$.

\begin{table}[t]
\caption{The following parameters were used in this example: $n=33, m=15, k=1,\ldots,14, 
\varepsilon =\frac{k}{n},$
  second column represents real probability of absolute deviation of the sample mean 
from $p = \frac{m}{n},$
  next two columns represent corresponding bound (\ref{eq:bnd56}) and Hoeffding's 
bound (\ref{eq:hff5}).}
   \label{tb:table1}
\begin{center}  \small
 \begin{tabular}{llll}
   \hline \noalign{\smallskip}
   $\varepsilon$ \hspace{0.3in} & True probability \hspace{0.02in} & Bound (\ref{eq:bnd56})
\hspace{0.05in} & Hoeffding's Bound (\ref{eq:hff5}) \\
\noalign{\smallskip}
\hline
\noalign{\smallskip}
0.0606  &       0.600713        &       0.78542 &       1.569446        \\
0.0909  &       0.382439        &       0.582175        &       1.159157        \\
0.1212  &       0.220522        &       0.38456 &       0.758396        \\
0.1515  &       0.114271        &       0.227376        &       0.43955 \\
0.1818  &       0.052796        &       0.120996        &       0.225672        \\
0.2121  &       0.021571        &       0.058319        &       0.102638        \\
0.2424  &       0.007724        &       0.025643        &       0.041352        \\
0.2727  &       0.0024  &       0.010366        &       0.014758        \\
0.3030  &       0.00064 &       0.003884        &       0.004666        \\
0.3333  &       0.000145        &       0.00136 &       \textbf{0.001307}        \\
0.3636  &       0.000027        &       0.000449        &       0.000324        \\
0.3939  &       0.000004        &       0.000141        &       0.000071        \\
0.4242  &       0.000001        &       0.000042        &       0.000014        \\
0.4545  &       0       &       0.000012        &       0.000002        \\
\hline
 \end{tabular}
\end{center}
\end{table}

\begin{remark} It follows from (\ref{eq:bnd63}) that the number of observations must
be big enough in order to ensure advantage of the Hoeffding's bound for the fixed
deviation parameter $\varepsilon$:
$$n > \frac{\log{2}}{2 \varepsilon^4}.$$
As a direct consequence, the bound (\ref{eq:bnd56}) will be so small:
$$  \exp{\{-\frac{2 \varepsilon^2 n}{1+\varepsilon^2}\}}
< \exp{\{-\frac{\log{2}}{\varepsilon^2(1+\varepsilon^2)}\}}$$
that any further improvement may not be regarded as a significant.
Figures~\ref{fig: figure1}(b)(d) illustrate above fact with relatively small numbers of
observations $n=20$ and $n=100.$
\end{remark}

\section{Continuous Case} \label{sec:asymptot}
\noindent
Suppose that $p$ and $\varepsilon$ are arbitrary numbers:
\begin{equation}
  \label{eq:ascond1}
    0< p < 1;
\end{equation}
\begin{equation}
  \label{eq:ascond2}
    \frac{1}{n} < \varepsilon \leq \min{\{p, 1-p\}}.
\end{equation}

\begin{remark} The special case $\varepsilon \leq \frac{1}{n}$ may be considered easily.
In the case $\varepsilon > p$ or $\varepsilon > 1-p$ we will have simplified cases
because we will need to approximate only one probability of deviation $\Pcal^+_2$ or 
$\Pcal^-_2$.
\end{remark}

We define a central point $np$ which is not necessarily integer.
We denote by 1) $h(n)$ - number of integer numbers in the left group
$[ np-n \varepsilon , np);$
2) $g(n)$ - number of integer numbers in the right group $[np, np+n\varepsilon].$
All remaining left groups will have $h$ integers with only one
possible exception as a last group.
Symmetrically, all remaining right groups will have $g$ integers with only one
possible exception as a last group.

Let us denote by $m$ the smallest integer in $[np, np+n\varepsilon].$

According to the construction
\begin{equation}
  \label{eq:ascond3}
    0 \leq \tilde{p} -p \leq \frac{1}{n}; \
\end{equation}
\begin{equation}
  \label{eq:ascond4}
    | \tilde{\varepsilon}_i - \varepsilon| \leq \frac{1}{n};
\end{equation}
\begin{equation}
  \label{eq:ascond5}
    \theta := \max{\{\frac{\varepsilon}{\tilde{\varepsilon}_1},
    \frac{\varepsilon}{\tilde{\varepsilon}_2}\}} 
\leq \frac{n \varepsilon}{n \varepsilon -1},
\end{equation}
where $\tilde{p} := \frac{m}{n}, \tilde{\varepsilon}_1 :=\frac{h(n)}{n},
 \tilde{\varepsilon}_2 :=\frac{g(n)}{n}.$

\begin{theorem} \label{th:teorem5} Suppose that conditions (\ref{eq:ascond1}) and
(\ref{eq:ascond2}) are valid. Then,
\begin{equation} \label{eq:funcbound}
   \Pb(|\frac{1}{n} \sum_{i=1}^n X_i - p| > \varepsilon)
   \leq \exp{\{ \varepsilon \cdot \varphi(n, \varepsilon)
   -\frac{2 \varepsilon^2 n}{1+\varepsilon^2}\}},
\end{equation}
where
\begin{equation} \label{eq:func42}
  \varphi(n, \varepsilon) :=
  \frac{\frac{n \varepsilon(1+\varepsilon^2)}{n \varepsilon-1}+
\varepsilon(2+6\varepsilon+\frac{9}{2n})}
  {(1+\varepsilon^2)(1+\varepsilon^2+\frac{1}{n}(1+3\varepsilon +\frac{9}{4n}))}.
\end{equation}
\end{theorem}
\textit{Proof:}
Again, we make an assumption $p \geq 0.5.$
Similar to (\ref{eq:bnd21}) we obtain
\begin{equation} \label{eq:asmp}
   \log{\lambda_1} \geq g \log{\frac{(1-p)(m+0.5g+1.5)}{p(n-m-0.5g-0.5)}}
 = g\log{(1+\delta)},
\end{equation}
where $\delta = \frac{m+1.5+0.5g-p-pn}{pn-mp-0.5p-0.5gp}.$

By (\ref{eq:gbound}), we obtain
\[ \frac{2(m+1.5+0.5g-p-pn)}{pn-2mp-gp-2p+m+1.5+0.5g} \geq \frac{\gamma g}{n} \]
or \hspace{0.5in}
$$ \tilde{\varepsilon}_2 +2(\tilde{p} -p) + \frac{3-2p}{n}
  \geq \gamma \tilde{\varepsilon}_2 \left( p - 2 \tilde{p} p -
  \tilde{\varepsilon}_2 p - \frac{2p}{n} + \tilde{p}
  + 0.5 \tilde{\varepsilon}_2  + \frac{3}{2n} \right). $$

Assuming that $\gamma \leq 2$
we simplify above inequality according to (\ref{eq:ascond3}) by making it stronger
$$\tilde{\varepsilon}_2 +\frac{3-2p}{n} \geq \gamma \tilde{\varepsilon}_2
\big( 2p(1-p) + \frac{1}{n}(1.5-2p) \big).$$
As a next step we can split above inequality into 2 inequalities
$$\tilde{\varepsilon}_2 \geq 2 \gamma \tilde{\varepsilon}_2 p(1-p); \hspace{0.05in}
\frac{3-2p}{n} \geq  \frac{\gamma \tilde{\varepsilon}_2}{n}(1.5-2p),$$
which are valid if
\begin{equation} \label{eq:asmp1}
  \gamma \leq 2.
\end{equation}

Now, we consider approximation in the left direction according to the previous framework
 (\ref{eq:bth63})
$$\log{q_1} \geq \frac{h}{2}
\log{\left( \left( \frac{p}{1-p}\right)^2 \frac{(n-m+2)(n-m+1+h)}{(m-h)(m-1)}\right)}$$
$$\geq \frac{h}{2} \log{\left(\frac{p(n-m+1+h)}{(1-p)(m-h)}\right)} = 
\frac{h}{2} \log{(1+\delta)},$$
where $\delta = \frac{pn+p-m+h}{m-mp-h+hp}.$

By (\ref{eq:gbound}), we obtain
$$\frac{pn+p-m+h}{m-2mp-h+2hp+pn+p} \geq \frac{\gamma h}{n} = 
\gamma \tilde{\varepsilon}_1$$
or
$$\tilde{\varepsilon}_1 +p-\tilde{p}+\frac{p}{n} \geq \gamma \tilde{\varepsilon}_1
 \left[ \tilde{p}(1-2p)+p- \tilde{\varepsilon}_1(1-2p) +\frac{p}{n}\right].$$

Again, we simplify above inequality  according to (\ref{eq:ascond3}) and (\ref{eq:ascond4})
\begin{equation} \label{eq:trnsf}
  \tilde{\varepsilon}_1 - \frac{1}{2n} \geq \gamma \tilde{\varepsilon}_1
  \left[ 2p(1-p) + \varepsilon (2p-1) + \frac{3p-1}{n}\right].
\end{equation}

Therefore,
\begin{equation} \label{eq:asmp2}
  \gamma \leq \frac{2 - \frac{\theta}{n \varepsilon}}{1+\varepsilon^2
  +\frac{1}{n}(1+3\varepsilon +\frac{9}{4n})}.
\end{equation}

The last remaining case corresponds directly to (\ref{eq:bnd31}).
$$\log{q_1} \geq h
\log{\left( \frac{p(n-m+1.5+0.5h)}{(1-p)(m-0.5h-0.5)}\right)} = h \log{(1+\delta)}$$
where $\delta=\frac{pn+p-m+0.5(h+1)}{m-pm+0.5(hp+p-h-1)}.$

By (\ref{eq:gbound}), we obtain \hspace{0.3in}
$$ \tilde{\varepsilon}_1 +2(p-\tilde{p})+\frac{2p+1}{n}
  \geq
   \gamma \tilde{\varepsilon}_1 \left[ \tilde{p}(1-2p) + 
p +\frac{\tilde{\varepsilon}_1}{2}(2p-1)
   + \frac{2p-0.5}{n}\right]. $$

We simplify above inequality  according to (\ref{eq:ascond3}) and (\ref{eq:ascond4})
\begin{equation} \label{eq:asmp3}
   \tilde{\varepsilon}_1 \geq \gamma  \tilde{\varepsilon}_1
  \left[ 2p(1-p) +\frac{\varepsilon (2p-1)}{2} + \frac{3p}{n} -\frac{1}{n}\right]
\end{equation}
which is valid for any
\begin{equation} \label{eq:bnd33a}
  0 < \gamma \leq 2 \left( \frac{4+\varepsilon^2}{4} +
  \frac{1}{n} \left(1+ \frac{3 \varepsilon}{2}+\frac{9}{4n}\right) \right)^{-1}.
\end{equation}

Finally, we derive required asymptotical relation as a consequence of the conditions
(\ref{eq:asmp1}), (\ref{eq:asmp2})
and (\ref{eq:bnd33a}) where condition (\ref{eq:asmp2}) is the most restrictive.
$\blacksquare$

\begin{table}[t]
\caption{Values of the function $\exp{\{\varepsilon \varphi(n, \varepsilon)\}}.$}
   \label{tb:table2}
\begin{center}  \small
 \begin{tabular}{lllllll}
   \hline \noalign{\smallskip}
$n \diagdown \varepsilon$ \hspace{0.1in} & 0.02 \hspace{0.01in} & 
0.05 \hspace{0.01in} &  0.1 \hspace{0.01in}
&  0.2 \hspace{0.01in} &  0.3 \hspace{0.01in} & 0.35     \\
\noalign{\smallskip}
\hline
\noalign{\smallskip}
$100$ & 1.0412  & 1.0697  &  1.1436  &  1.3737  & 1.7612  & 2.0357  \\
$1,000$  & 1.0221  & 1.0582  & 1.1336  &  1.3652  &       1.7566  &       2.0349  \\
$100,000$  & 1.0211  & 1.0572  & 1.1326  & 1.3643  & 1.7561  &  2.0348  \\
\hline
 \end{tabular}
\end{center}
\end{table}

\subsection{Asymptotical Bounds for Normalized Sum of Bernoulli 
Random Variables} \label{ssec:normal}
\noindent
Let us denote by $F_n$ distribution function of the normalized 
sum of Bernoulli random variables
\begin{equation} \label{eq:norm}
  \eta_n = \frac{\sum_{i=1}^n X_i -p}{\sqrt{np(1-p)}}.
\end{equation}

According to the Central Limit Theorem,
$$\lim_{n \rightarrow \infty} F_n(t) = \Phi(t),$$
where $\Phi(t)$ is a  standard normal distribution function.
Respectively, it appears to be reasonable to use random variable (\ref{eq:norm}) as a test case
in order to compare different bounds.

\begin{proposition} \label{pr:prop3} The bounds (\ref{eq:bnd1}), (\ref{eq:bnd56}) and
(\ref{eq:funcbound}), as an upper bounds for the following probability $\Pb(|\eta_n| > t),$
are equivalent asymptotically and equal to
\begin{equation} \label{eq:clt7}
 \alpha \exp{\{-2 t^2 p(1-p)\}},
\end{equation}
where $\alpha =1,$ and $0 \leq t, p \leq 1.$
\end{proposition}
\textit{Proof:} 
We have
$$\Pb(|\eta_n| > t) = \Pb \left( |\frac{1}{n} \sum_{i=1}^n X_i -p| >
  t_n \right),$$
where $t_n = t \sqrt{\frac{p(1-p)}{n}} $.

Then, we insert $t_n$ into (\ref{eq:bnd1}), (\ref{eq:bnd56}) or
(\ref{eq:funcbound}), and take $lim$ if $n \rightarrow \infty.$
As a result, we obtain required formulae (\ref{eq:clt7}).
$\blacksquare$

\begin{remark} Using Hoeffding's bound we will obtain the same asymptotical structure
(\ref{eq:clt7}), but with  $\alpha=2.$
\end{remark}

\section{Concluding remarks}  \label{sec:final}
The Proposition~\ref{pr:prop} proves that value of
the linear coefficient $\alpha=1$ can not be improved. Taking this fact as a base point we
established a new structure of the exponential bound (\ref{eq:bnd56}).
Figures~\ref{fig: figure1}(a),(c), (d) and Table~\ref{tb:table1}
demonstrate advantages of the bound (\ref{eq:bnd56})
against Hoeffding's bound if value of $\varepsilon$ is small enough.
The above Section~\ref{ssec:normal} demonstrates additional arguments in support of
the proposed bounds, and these arguments cover not to only Bernoulli
random variables.
The area of applications of the proposed methods may be extended further, for example,
we can consider arbitrary bounded random variables or uniform metric for the
empirical distributions.

The paper represents a fresh look at the ideas and methods which
Ja. Bernoulli proposed in the 17th century. As it was demonstrated
in the Section~\ref{sec:bernoulli} the probability of deviation in
the classical Bernoulli case may be bounded using standard
structure of the exponential bound with optimal linear coefficient
$\alpha=1$. It was the first step of this research which was
completed in 1987 shortly after the 1st World Congress of the
Bernoulli Society where the author purchased book  \cite{Bern86}.

\end{document}